\input amstex
\input amsppt.sty
 \magnification=\magstep1 \hsize=36truecc
 \vsize=23.5truecm
\baselineskip=14truept
\def\q{\quad}
\def\qq{\qquad}
\def\mod#1{\ (\text{\rm mod}\ #1)}

\def\t{\text}
\def\qtq#1{\q\t{#1}\q}
\par The paper will appear in Proc. Amer. Math.
Soc..\par\q
\def\mod#1{\ (\text{\rm mod}\ #1)}
\def\qtq#1{\q\t{#1}\q}
\def\f{\frac}
\def\e{\equiv}
\def\b{\binom}

\def\sls#1#2{(\f{#1}{#2})}
 \def\ls#1#2{\big(\f{#1}{#2}\big)}
\def\Ls#1#2{\Big(\f{#1}{#2}\Big)}
\let \pro=\proclaim
\let \endpro=\endproclaim
\topmatter
\title Congruences concerning Legendre polynomials
\endtitle
\author ZHI-Hong Sun\endauthor
\affil School of Mathematical Sciences, Huaiyin Normal University,
\\ Huaian, Jiangsu 223001, PR China
\\ E-mail: szh6174$\@$yahoo.com
\\ Homepage: http://www.hytc.edu.cn/xsjl/szh
\endaffil

 \nologo \NoRunningHeads

\abstract{Let $p$ be an odd prime. In the paper, by using the
properties of Legendre polynomials we prove some congruences for
$\sum_{k=0}^{\f{p-1}2}\b{2k}k^2m^{-k}\mod {p^2}$. In particular, we
confirm several conjectures of Z.W. Sun. We also pose $13$
conjectures on supercongruences.
\par\q
\newline MSC: Primary 11A07, Secondary 33C45, 11E25 \newline Keywords:
Legendre polynomial; Congruence}
 \endabstract
  \footnote"" {The author is
supported by the Natural Sciences Foundation of China (grant
10971078).}
\endtopmatter
\document
\subheading{1. Introduction}
\par\par\q Let $p$ be an odd prime. In 2003, Rodriguez-Villegas [11]
 conjectured the following congruence:
$$\sum_{k=0}^{(p-1)/2}\f{\b{2k}k^2}{16^k}\e
(-1)^{\f{p-1}2}\mod{p^2}.\tag 1.1$$ This was later confirmed by
Mortenson [7] via the Gross-Koblitz formula. See also [9] and [10,
p.204]. Recently my twin brother Zhi-Wei Sun [13] obtained the
congruences for $\sum_{k=0}^{p-1}\b{2k}k^2m^{-k}\mod p$ in the cases
$m=8,-16,32$, and made several conjectures for
$\sum_{k=0}^{p-1}\b{2k}k^2m^{-k}\mod {p^2}$. For example, he
conjectured $$ \sum_{k=0}^{(p-1)/2}\f{\b{2k}k^2}{32^k} \e\cases
0\mod {p^2}&\t{if $4\mid p-3$,}
\\2a-\f p{2a}\mod {p^2}&\t{if $4\mid p-1$ and $p=a^2+b^2$ with
$4\mid a-1$.}
\endcases\tag 1.2$$
\par Let $\{P_n(x)\}$ be the Legendre polynomials given by
$$\f 1{\sqrt{1-2xt+t^2}}=\sum_{n=0}^{\infty}P_n(x)t^n\q(|t|<1).$$ It
is well known that (see [6, pp.\;228-232], [4, (3.132)-(3.133)])
$$P_n(x)=\f
1{2^n}\sum_{k=0}^{[n/2]}\f{(-1)^k(2n-2k)!}{k!(n-k)!(n-2k)!}x^{n-2k}
=\f 1{2^n\cdot n!}\cdot\f{d^n}{dx^n}(x^2-1)^n\tag 1.3$$ and
$(n+1)P_{n+1}(x)=(2n+1)xP_n(x)-nP_{n-1}(x),$ where $[x]$ is the
greatest integer not exceeding $x$.
\par In the paper, by using the expansions of Legendre polynomials
we  obtain some congruences for $P_{\f{p-1}2}(x)$ modulo $p^2$,
where $p$ is an odd prime and $x$ is a rational $p$-integer. For
example, we have
$$\sum_{k=0}^{p-1}\f{\b{2k}k^2}{16^k}\big(x^k-(-1)^{\f{p-1}2}(1-x)^k\big)
\e 0\mod{p^2},\tag 1.4$$ and
$$\sum_{k=0}^{\f{p-1}2}\f{\b{2k}k^2}{16^k}\Big(x^k-\Ls
{x}px^{-k}\Big)\e 0\mod p\qtq{for}x\not\e 0\mod p,\tag 1.5$$ where
$\ls xp$ is the Legendre symbol.  Taking $x=1$ in (1.4) we obtain
(1.1) immediately, and taking $x=\f 12$ in (1.4) we deduce (1.2) for
$p\e 3\mod 4$. We also determine
$\sum_{k=0}^{\f{p-1}2}\f{\b{2k}k^2}{32^k}k(k-1)\cdots(k-r+1)\mod
{p^2}$ for $r\in\{1,2,\ldots,\f{p-1}2\}$,
$\sum_{k=0}^{[p/3]}\f{(3k)!}{54^k\cdot k!^3}\mod p$ and pose some
conjectures on supercongruences concerning binary quadratic forms.
\par Throughout this paper we use $\Bbb Z,\ \Bbb N$ and $\Bbb Z_p$ to denote the sets of
integers, positive integers and rational $p$-integers for a prime
$p$, respectively.
\subheading{2. Main results}
\par\q  \pro{Lemma 2.1} For $n\in\Bbb N$ we have
$$P_n(x)=\sum_{k=0}^n\b{n+k}{2k}\b {2k}k\Ls{x-1}2^k.$$
\endpro
Proof. From [4, (3.135)] we have the following result due to Murphy:
$$P_n(x)=\sum_{k=0}^n\b nk\b {n+k}k\Ls{x-1}2^k.\tag 2.1$$
As $\b{n+k}{2k}\b {2k}k=\b {n+k}k\b nk$, we obtain the result.

\pro{Lemma 2.2} Let $p$ be an odd prime and
$k\in\{1,2,\ldots,(p-1)/2\}$. Then
$$\b{\f{p-1}2+k}{2k}\e \f{\b{2k}k}{(-16)^k}\Big(1-p^2
\sum_{i=1}^k\f 1{(2i-1)^2}\Big)\mod {p^4}.$$
\endpro
Proof. Clearly
$$\aligned
\b{\f{p-1}2+k}{2k}&=\f{(\f{p-1}2+k)(\f{p-1}2+k-1)\cdots(\f{p-1}2-k+1)}{(2k)!}
\\&=\f{(p+2k-1)(p+2k-3)\cdots(p-(2k-3))(p-(2k-1))}{2^{2k}\cdot
(2k)!}
\\&=\f{(p^2-1^2)(p^2-3^2)\cdots (p^2-(2k-1)^2)}{2^{2k}\cdot
(2k)!}
\\&\e \f {(-1)^k\cdot 1^2\cdot 3^2\cdots
(2k-1)^2}{2^{2k}\cdot (2k)!}\Big(1-p^2\sum_{i=1}^k\f
1{(2i-1)^2}\Big)\mod {p^4}.\endaligned$$ To see the result, we note
that
$$\f{1^2\cdot 3^2\cdots(2k-1)^2}{2^{2k}\cdot
(2k)!}=\f{(2k)!^2}{(2\cdot  4\cdots (2k))^2\cdot 2^{2k}\cdot
(2k)!}=\f{(2k)!}{2^{4k}\cdot k!^2}=\f{\b{2k}k}{16^k}.$$

\par Let $p$ be an odd prime, and let $\{A(n)\}$ be the Ap\'ery numbers given by
$$A(n)=\sum_{k=0}^n\b{n+k}k^2\b nk^2.$$
It is well known that (see [1],[10]) $A(\f{p-1}2)\e a(p)\mod {p^2},$
where $a(n)$ is defined by
$$q\prod_{n=1}^{\infty}(1-q^{2n})^4(1-q^{4n})^4=\sum_{n=1}^{\infty}a(n)q^n.$$
 By the fact $\b{n+k}k\b nk=\b{n+k}{2k}\b{2k}k$ and Lemma 2.2 we have
$$A\Ls{p-1}2=\sum_{k=0}^{\f{p-1}2}\b{\f{p-1}2+k}{2k}^2\b{2k}k^2
\e
\sum_{k=0}^{\f{p-1}2}\Big(\f{\b{2k}k}{(-16)^k}\Big)^2\b{2k}k^2\mod
{p^2}.$$ Hence
$$a(p)\e A\Ls{p-1}2\e
\sum_{k=0}^{\f{p-1}2}\f{\b{2k}k^4}{4^{4k}}\mod{p^2}.\tag 2.2$$ Let
$b(n)$ be given by
$q\prod_{n=1}^{\infty}(1-q^{4n})^6=\sum_{n=1}^{\infty}b(n)q^n.$ Then
Mortenson [8] proved the following conjecture of Rodriguez-Villegas:
$$\sum_{k=1}^{\f{p-1}2}\f{\b{2k}k^3}{4^{3k}}\e b(p)\mod{p^2}.\tag 2.3$$

\pro{Theorem 2.1} Let $p$ be an odd prime and let $x$ be a variable.
Then
$$\sum_{k=0}^{p-1}\f{\b{2k}k^2}{16^k}\big(x^k-(-1)^{\f{p-1}2}(1-x)^k\big)
\e
\sum_{k=0}^{\f{p-1}2}\f{\b{2k}k^2}{16^k}\big(x^k-(-1)^{\f{p-1}2}(1-x)^k\big)\e
0\mod {p^2}.$$\endpro Proof. For a variable $t$, by Lemmas 2.1 and
2.2 we have
$$ P_{\f{p-1}2}(t)=\sum_{k=0}^{\f{p-1}2}\b{\f{p-1}2+k}{2k}\b {2k}k\Ls{t-1}2^k
\e \sum_{k=0}^{\f{p-1}2}\f{\b{2k}k^2}{16^k}\Ls{1-t}2^k\mod
{p^2}.\tag 2.4$$ It is known that (see [6]) $P_n(t)=(-1)^nP_n(-t)$.
Thus, by (2.4),
$$\sum_{k=0}^{\f{p-1}2}\f{\b{2k}k^2}{16^k}\Ls{1-t}2^k\e
(-1)^{\f{p-1}2}\sum_{k=0}^{\f{p-1}2}\f{\b{2k}k^2}{16^k}\Ls{1+t}2^k\mod{p^2}.$$
Now taking $t=1-2x$ in the congruence we deduce
$\sum_{k=0}^{\f{p-1}2}\f{\b{2k}k^2}{16^k}\big(x^k-(-1)^{\f{p-1}2}(1-x)^k\big)\e
0\mod {p^2}.$ To complete the proof, we note that for
$k\in\{\f{p+1}2,\f{p+3}2,\ldots,p-1\}$,
$\b{2k}k=2k(2k-1)\cdots(k+1)/k!\e 0\mod p$.
 \pro{Theorem 2.2} Let $p$ be an odd prime. Then
$$ \sum_{k=0}^{\f{p-1}2}\f{\b{2k}k^2}{32^k} \e\cases 0\mod
{p^2}&\t{if $4\mid p-3$,}
\\2a-\f p{2a}\mod {p^2}&\t{if $4\mid p-1$ and $p=a^2+b^2$ with
$4\mid a-1$.}
\endcases$$
\endpro
Proof. When $p\e 3\mod 4$, taking $x=\f 12$ in Theorem 2.1 we obtain
the result. Now suppose $p\e 1\mod 4$ and so $p=a^2+b^2$ with
$a,b\in\Bbb Z$ and $a\e 1\mod 4$. It is well known that ([6])
$$P_{2n+1}(0)=0\qtq{and}P_{2n}(0)=\f{(-1)^n}{2^{2n}}\b{2n}n.\tag 2.5$$
Thus, by (2.4) and (2.5) we have
$$\sum_{k=0}^{\f{p-1}2}\f{\b{2k}k^2}{32^k}\e
P_{\f{p-1}2}(0)=\f{(-1)^{\f{p-1}4}}{2^{\f{p-1}2}}\b{\f{p-1}2}{\f{p-1}4}\mod
{p^2}.$$ According to the result due to Chowla, Dwork and Evans (see
[2] or [3]), we have
$$\b{\f{p-1}2}{\f{p-1}4}\e
\f{2^{p-1}+1}2\Big (2a-\f p{2a}\Big)\mod {p^2}.$$ Set
$q=(2^{\f{p-1}2}-(-1)^{\f{p-1}4})/p$. Then $2^{p-1}\e
1+2(-1)^{\f{p-1}4}qp\mod {p^2}$. Thus
$$\f{2^{p-1}+1}{2\cdot 2^{\f{p-1}2}}\e
\f{2+2(-1)^{\f{p-1}4}qp}{2((-1)^{\f{p-1}4}+qp)}=(-1)^{\f{p-1}4}\mod{p^2}.$$
Hence
$$\sum_{k=0}^{\f{p-1}2}\f{\b{2k}k^2}{32^k}\e
\f{(-1)^{\f{p-1}4}}{2^{\f{p-1}2}}\b{\f{p-1}2}{\f{p-1}4} \e
\f{(-1)^{\f{p-1}4}}{2^{\f{p-1}2}}\cdot \f{2^{p-1}+1}2\Big (2a-\f
p{2a}\Big)\e 2a-\f p{2a}\mod {p^2}.$$ The proof is now complete.
 \par{\bf Remark 2.1} Theorem 2.2 was conjectured by Zhi-Wei Sun
 ([13]),  and the congruence for
  $\sum_{k=0}^{\f{p-1}2}\f{\b{2k}k^2}{32^k}\mod p$ was also proved by
  Zhi-Wei Sun in [13].

 \pro{Theorem 2.3} Let $p$ be an odd prime and $r\in\{1,2,\ldots,(p-1)/2\}$. Then
 $$\aligned&\sum_{k=0}^{\f{p-1}2}\f{\b{2k}k^2}{32^k}k(k-1)\cdots(k-r+1)
\\&\e \cases 0\mod {p^2}&\t{if $4\mid (p+1-2r)$,}
\\(-1)^{\f{p-1+2r}4}2^{-\f{p-1}2}\f{(\f{p-1}2+r)!}{\f{p-1-2r}4!\f{p-1+2r}4!}
\mod{p^2}&\t{if $4\mid (p-1-2r)$.}
\endcases\endaligned$$
\endpro
Proof. By (2.4) we have
$$\aligned\f{d^r\;P_{\f{p-1}2}(t)}{dt^r}&\e\sum_{k=0}^{\f{p-1}2}\f{\b{2k}k^2}{(-32)^k}\cdot
\f{d^r(t-1)^k}{dt^r}
\\&=\sum_{k=0}^{\f{p-1}2}\f{\b{2k}k^2}{(-32)^k}k(k-1)\cdots(k-r+1)(t-1)^{k-r}
\mod{p^2}.\endaligned\tag 2.6$$ Hence
$$\f{d^r\;P_{\f{p-1}2}(t)}{dt^r}\bigg|_{t=0}=(-1)^r\sum_{k=0}^{\f{p-1}2}\f{\b{2k}k^2}{32^k}
k(k-1)\cdots(k-r+1).$$ By (1.3) we have
$$\aligned\f{d^r}{dt^r}P_{\f{p-1}2}(t)&=\f
1{2^{(p-1)/2}}\cdot\f{d^r}{dt^r}
\sum_{m=0}^{[\f{p-1}4]}\f{(-1)^m(p-1-2m)!}{m!(\f{p-1}2-m)!(\f{p-1}2-2m)!}t^{\f{p-1}2-2m}
\\&=\f 1{2^{(p-1)/2}}
\sum_{m=0}^{[\f{p-1-2r}4]}\f{(-1)^m(p-1-2m)!}{m!(\f{p-1}2-m)!(\f{p-1}2-2m)!}
\\&\qq\times(\f{p-1}2-2m)(\f{p-1}2-2m-1)\cdots(\f{p-1}2-2m-r+1)t^{\f{p-1}2-2m-r}.
\endaligned$$
Thus,
$$\f{d^rP_{\f{p-1}2}(t)}{dt^r}\Big|_{t=0}=\cases 0&\t{if $r\not\e
\f{p-1}2\mod 2$,}
\\\f{(-1)^m(p-1-2m)!}{2^{(p-1)/2}\cdot m!(\f{p-1}2-m)!}&\t{if
$r=\f{p-1}2-2m$.}\endcases$$ Now combining all the above we obtain
the result.
 \pro{Corollary 2.1} Let $p$ be an odd prime. Then
$$\sum_{k=0}^{\f{p-1}2}\f{k^2\b{2k}k^2}{32^k}
\e\cases
(-1)^{\f{p+3}4}2^{-\f{p-1}2}\f{\f{p+3}2!}{\f{p-5}4!\f{p+3}4!}\mod{p^2}
&\t{if $p\e 1\mod 4$,}
\\(-1)^{\f{p+1}4}2^{-\f{p-1}2}\f{\f{p+1}2!}{\f{p-3}4!\f{p+1}4!}\mod{p^2}
&\t{if $p\e 3\mod 4$.}
\endcases$$
\endpro
Proof. By Theorem 2.3 we have
$$\sum_{k=0}^{\f{p-1}2}\f{k\b{2k}k^2}{32^k}
\e\cases 0\mod{p^2} &\t{if $p\e 1\mod 4$,}
\\(-1)^{\f{p+1}4}2^{-\f{p-1}2}\f{\f{p+1}2!}{\f{p-3}4!\f{p+1}4!}\mod{p^2}
&\t{if $p\e 3\mod 4$}
\endcases$$
and
$$\sum_{k=0}^{\f{p-1}2}\f{k(k-1)\b{2k}k^2}{32^k}
\e\cases
(-1)^{\f{p+3}4}2^{-\f{p-1}2}\f{\f{p+3}2!}{\f{p-5}4!\f{p+3}4!}\mod{p^2}
&\t{if $p\e 1\mod 4$,}
\\0\mod{p^2}
&\t{if $p\e 3\mod 4$.}
\endcases$$
Observe that $k^2=k(k-1)+k$. From the above we deduce the result.

 \pro{Lemma 2.3} Let $p$ be a prime greater than $3$ and let $t$ be a variable. Then
 $$P_{[\f p3]}(t)\e \sum_{k=0}^{[p/3]}\f{(3k)!}{k!^3}\Ls{1-t}{54}^k\mod
 p.$$
 \endpro
 Proof. Suppose $r=1$ or $2$ according as $3\mid
p-1$ or $3\mid p-2$. Then clearly
$$\aligned \b{\f {p-r}3+k}{2k}&
=\f{(\f{p-r}3+k)(\f{p-r}3+k-1)\cdots(\f{p-r}3-k+1)}{(2k)!}
\\&=\f{(p+3k-r)(p+3k-r-3)\cdots(p-(3k+r-3))}{3^{2k}\cdot (2k)!}
\\&\e (-1)^k\f{(3k-r)(3k-r-3)\cdots(3-r)\cdot r(r+3)\cdots(3k+r-3)}
{3^{2k}\cdot (2k)!}
\\&=\f{(-1)^k\cdot (3k)!}{3\cdot 6\cdots 3k\cdot 3^{2k}\cdot (2k)!}
=\f{(-1)^k\cdot(3k)!}{3^k\cdot k!\cdot 3^{2k}\cdot (2k)!} \mod p.
\endaligned$$
Hence, by Lemma 2.1 we have
$$\aligned P_{[\f p3]}(t)&=\sum_{k=0}^{[p/3]}\b{[\f
p3]+k}{2k}\b{2k}k\ls{t-1}2^k  \e \sum_{k=0}^{[p/3]}\f{(-1)^k\cdot
(3k)!}{3^{3k}\cdot k!(2k)!}\cdot \f{(2k)!}{k!^2}\Ls{t-1}2^k
\\&= \sum_{k=0}^{[p/3]}\f{(3k)!}{27^k\cdot k!^3}\Ls{1-t}2^k\mod p.\endaligned$$
This proves the lemma.
 \pro{Theorem
2.4} Let $p$ be a prime greater than $3$ and let $x$ be a variable.
Then
$$\sum_{k=0}^{[p/3]}\f{(3k)!}{27^k\cdot k!^3}\big(x^k-(-1)^{[\f p3]}(1-x)^k\big)
\e 0\mod p.$$\endpro Proof. As $P_n(t)=(-1)^nP_n(-t)$, using Lemma
2.3 we deduce
$$\sum_{k=0}^{[p/3]}\f{(3k)!}{27^k\cdot k!^3}\Big(\Ls{1-t}2^k-(-1)^{[p/3]}\Ls{1+t}2^k
\Big)\e 0\mod p.$$ Now putting $t=1-2x$ in the congruence we obtain
the result. \pro{Corollary 2.2} Let $p$ be a prime greater than $3$.
Then
$$\sum_{k=0}^{[p/3]}\f{(3k)!}{27^k\cdot k!^3}\e \Ls p3\mod p.$$
\endpro
Proof. Taking $x=1$ in Theorem 2.4 and noting that
$(-1)^{[p/3]}=\sls p3$ we deduce the result.
\par{\bf Remark 2.2} By [8] or [10, p.\;204] we have the following
stronger supercongruence $\sum_{k=0}^{p-1}\f{(3k)!}{27^k\cdot
k!^3}\e \sls p3\mod {p^2}.$
 \pro{Lemma 2.4} Let $p$
be an odd prime and $k\in\{1,2,\ldots,\f{p-1}2\}$. Then
$$\b{(p-1)/2}k\e \f 1{(-4)^k}\b{2k}k\Big(1-p\sum_{i=1}^k\f 1{2i-1}\Big)\mod {p^2}.$$
\endpro
Proof. It is clear that
$$\aligned\b{\f{p-1}2}k&=\f{\f{p-1}2(\f{p-1}2-1)\cdots(\f{p-1}2-k+1)}{k!}
=\f{(p-1)(p-3)\cdots(p-(2k-1))}{2^k\cdot k!}
\\&\e \f{(-1)(-3)\cdots(-(2k-1))}{2^k\cdot k!}\Big(1-p\sum_{i=1}^k\f 1{2i-1}\Big)
\\&=\f{(-1)^k\cdot (2k)!}{(2^k\cdot k!)^2}\Big(1-p\sum_{i=1}^k\f 1{2i-1}\Big)
\mod {p^2}.\endaligned$$ This yields the result.

 \pro{Theorem 2.5} Let $p$ be a prime greater
than $5$. Then
$$\sum_{k=0}^{[p/3]}\f{(3k)!}{54^k\cdot k!^3}\e \cases 0\mod p&\t{if $6\mid p-5$,}
\\2A\mod p&\t{if $6\mid p-1$ and $p=A^2+3B^2$ with $3\mid A-1$}
\endcases$$
and
$$\sum_{k=0}^{[p/3]}\f{k\cdot (3k)!}{54^k\cdot k!^3}\e \cases 0\mod p&\t{if $6\mid p-1$,}
\\\f 13(-1)^{\f{p+1}6} 2^{-\f{p+1}3}\b{(p+1)/3}{(p+1)/6}\mod p&\t{if $6\mid p-5$.}
\endcases$$
\endpro
Proof. Taking $t=0$ in Lemma 2.3 and applying (2.5) and Lemma 2.4 we
deduce that
$$\sum_{k=0}^{[p/3]}\f{(3k)!}{54^k\cdot k!^3}
\e \cases 0\mod p&\t{if $p\e 5\mod
6$,}\\\f{(-1)^{(p-1)/6}}{2^{(p-1)/3}}\b{(p-1)/3}{(p-1)/6}\e
\b{(p-1)/2}{(p-1)/6}\mod p&\t{if $p\e 1\mod 6$.}\endcases$$ Now
suppose $p\e 1\mod 6$ and so $p=A^2+3B^2$ with $A,B\in\Bbb Z$ and
$A\e 1\mod 3$. By [2, Theorem 9.4.4] we have $\b{(p-1)/2}{(p-1)/6}\e
2A\mod p.$ Thus the first part follows.
\par By Lemma 2.3 we have
$$\f {d}{dt}P_{[\f p3]}(t)\e -\sum_{k=0}^{[p/3]}\f{(3k)!}{54^k\cdot k!^3}\cdot
k(1-t)^{k-1}\mod p.$$ Thus, $\f {d}{dt}P_{[\f p3]}(t)\big|_{t=0}\e
-\sum_{k=0}^{[p/3]}\f{k\cdot (3k)!}{54^k\cdot k!^3}\mod p.$ From
(1.3) we know that
$$\f {d}{dt}P_{[\f p3]}(t)\big|_{t=0}=\cases 0&\t{if $p\e 1\mod 6$,}
\\ 2^{-\f{p-2}3}\cdot(-1)^{\f{p-5}6}\f{\f{p+1}3!}{\f{p-5}6!\f{p+1}6!}
&\t{if $p\e 5\mod 6$.}
\endcases$$
Thus the second part is true.

 \pro{Lemma 2.5} Let $p$ be an odd prime and
$k\in\{1,2,\ldots,\f{p-1}2\}$. Then
$$\f{(-1)^k\b{(p-1)/2+k}k}{\b{(p-1)/2}k}\e 1+2p\sum_{i=1}^k\f 1{2i-1}
\e 3-2(-4)^k\f{\b{(p-1)/2}k}{\b{2k}k}\mod {p^2}.$$
\endpro
Proof. It is clear that
$$\aligned \f{(-1)^k\b{(p-1)/2+k}k}{\b{(p-1)/2}k}&
=\f{(\f{p-1}2+k)(\f{p-1}2+k-1)\cdots(\f{p-1}2+1)}{(-1)^k\f{p-1}2(\f{p-1}2-1)\cdots(\f{p-1}2-k+1)}
\\&=\f{(p+2k-1)(p+2k-3)\cdots(p+1)}{(-1)^k(p-1)(p-3)\cdots(p-(2k-1))}
\\&\e \f{1\cdot 3\cdots(2k-1)(1+p\sum_{i=1}^k\f 1{2i-1})}
{1\cdot 3\cdots(2k-1)(1-p\sum_{i=1}^k\f 1{2i-1})}
\\&\e \Big(1+p\sum_{i=1}^k\f 1{2i-1}\Big)^2\e 1+2p\sum_{i=1}^k\f 1{2i-1}
\mod {p^2}.\endaligned$$ This together with Lemma 2.4 yields the
result. \pro{Theorem 2.6} Let $p$ be an odd prime, $x\in\Bbb Z_p$
and $x\not\e -1\mod p$. Then
$$P_{\f{p-1}2}(x)\e \Ls{2(x+1)}pP_{\f{p-1}2}\Ls{3-x}{1+x}\mod p.$$
\endpro
Proof. It is known that (see [4, (3.134)])
$$P_n(x)=\sum_{k=0}^n\b nk^2\Ls{x+1}2^{n-k}\Ls{x-1}2^k.$$
Thus, using Lemma 2.5 and (2.1) we see that
$$\aligned P_{\f{p-1}2}(x)&=\Ls{x+1}2^{\f{p-1}2}
\sum_{k=0}^{\f{p-1}2}\b{\f{p-1}2}k^2\Ls{x-1}{x+1}^k
\\&\e \Ls {2(x+1)}p\sum_{k=0}^{\f{p-1}2}\b{\f{p-1}2}k\b{\f{p-1}2+k}k(-1)^k\Ls{x-1}{x+1}^k
\\&=\Ls {2(x+1)}pP_{\f{p-1}2}\Big(1+2\cdot \f{1-x}{1+x}\Big)=\Ls {2(x+1)}p
P_{\f{p-1}2}\Ls{3-x}{1+x}\mod p.\endaligned$$ This proves the
theorem.
 \pro{Corollary 2.3} Let $p$ be a prime of the form $4k+3$.
Then $p\mid P_{\f{p-1}2}(3)$.
\endpro
Proof. By Theorem 2.6 and (2.5) we have $P_{\f{p-1}2}(3)\e \sls
2pP_{\f{p-1}2}(0)=0\mod p$.

 \pro{Theorem 2.7} Let $p$ be an odd prime, $x\in\Bbb Z_p$ and
$x\not\e 0\mod p$. Then
$$\sum_{k=0}^{\f{p-1}2}\f{\b{2k}k^2}{16^k}\Big(x^k-\Ls
{x}px^{-k}\Big)\e 0\mod p.$$
\endpro
Proof. Clearly the result is true for $x\e 1\mod p$.  Now assume
$x\not\e 1\mod p$. As $P_n(t)=(-1)^nP_n(-t)$ (see [6]), using
Theorem 2.6 we see that for $t\in\Bbb Z_p$ with $t\not\e \pm 1\mod
p$,
$$\ls{2(t+1)}pP_{\f{p-1}2}\Ls{3-t}{1+t}\e
(-1)^{\f{p-1}2}\Ls{2(-t+1)}pP_{\f{p-1}2}\Ls{3+t}{1-t}\mod p.$$ Thus,
$$P_{\f{p-1}2}\Ls{3-t}{1+t}\e \Ls{t^2-1}pP_{\f{p-1}2}\Ls{3+t}{1-t}\mod p.$$
Now applying (2.4) we deduce that
$$\sum_{k=0}^{\f{p-1}2}\f{\b{2k}k^2}{16^k}\Ls{1-\f{3-t}{1+t}}2^k
\e\Ls{t^2-1}p\sum_{k=0}^{\f{p-1}2}\f{\b{2k}k^2}{16^k}\Ls{1-\f{3+t}{1-t}}2^k
\mod p$$ and so
$$\sum_{k=0}^{\f{p-1}2}\f{\b{2k}k^2}{16^k}\Big(\Ls{t-1}{t+1}^k-\Ls{(t-1)/(t+1)}p\Ls{t+1}{t-1}^k
\Big)\e 0\mod p.$$ Set $t=(1+x)/(1-x)$. Then $t\not\e \pm 1\mod p$
and $x=(t-1)/(t+1)$. Hence the result follows.
\par\q
\par Let $p$ be an odd prime, and $x\in\Bbb Z_p$ with $x\not\e 0,1\mod p$. By Theorems
2.1 and 2.7 we have
$$\sum_{k=0}^{(p-1)/2}\f{\b{2k}k^2}{(16x)^k}
\e\Ls{x(x-1)}p\sum_{k=0}^{(p-1)/2}\f{\b{2k}k^2}{(16(1-x))^k}\mod
p.\tag 2.7$$

 \pro{Theorem 2.8} Let $p$ be an odd prime, $x\in\Bbb
Z_p$ and $x\not\e 0\mod p$. Then
$$\sum_{k=0}^{\f{p-1}2}\b{2k}k^2\Ls x4^{2k}\e \Ls{-x}pP_{\f{p-1}2}
\Ls{x+x^{-1}}2\mod p.$$
\endpro
Proof. From [4, (3.138)] we have the following result due to
Kelisky:
$$\sum_{k=0}^n\b{2n-2k}{n-k}\b{2k}kx^{2k}=2^{2n}x^nP_n\Ls{x+x^{-1}}2.\tag 2.8$$
Taking $n=(p-1)/2$ in (2.8) we have
$$\sum_{k=0}^{\f{p-1}2}\b{p-1-2k}{\f{p-1}2-k}\b{2k}kx^{2k}
=2^{p-1}x^{\f{p-1}2}P_{\f{p-1}2}\Ls{x+x^{-1}}2\e \Ls
xpP_{\f{p-1}2}\Ls{x+x^{-1}}2\mod p.$$ To see the result, using Lemma
2.2 we note that for $0\le k\le \f{p-1}2$,
$$\align\b{p-1-2k}{\f{p-1}2-k}&=\f{(p-1-2k)(p-2-2k)\cdots
(p-(\f{p-1}2+k))}{(\f{p-1}2-k)!}\tag 2.9\\& \e
(-1)^{\f{p-1}2-k}\f{(2k+1)(2k+2)\cdots (\f{p-1}2+k)}{(\f{p-1}2-k)!}
\\&=(-1)^{\f{p-1}2-k}\b{\f{p-1}2+k}{2k}\e
(-1)^{\f{p-1}2}\f{\b{2k}k}{16^k}\mod p.\endalign$$

 \pro{Theorem 2.9} Let $p$ be a prime of the form $4k+1$ and $p=a^2+b^2$ with $a,b\in\Bbb Z$ and $a\e 1\mod 4$.
  Then
 $$\sum_{k=0}^{\f{p-1}2}\f{\b{2k}k^2}{8^k}\e
 \sum_{k=0}^{\f{p-1}2}\f{\b{2k}k^2}{(-16)^k}
 \e P_{\f{p-1}2}(3)\e
 (-1)^{\f{p-1}4}\Big(2a-\f p{2a}\Big)\mod {p^2}.$$
 \endpro
 Proof. By Theorem 2.1 and (2.4) we have
$\sum_{k=0}^{\f{p-1}2}\b{2k}k^28^{-k}\e
 \sum_{k=0}^{\f{p-1}2}\b{2k}k^2(-16)^{-k}\e P_{\f{p-1}2}(3)$ $\mod {p^2}.$
From Theorem 2.6, (2.5) and Gauss' congruence
$\b{(p-1)/2}{(p-1)/4}\e 2a$ $\mod p$ (see [2]) we have
$$P_{\f{p-1}2}(3)\e \ls 2pP_{\f{p-1}2}(0)=\Ls
2p\f{(-1)^{(p-1)/4}}{2^{(p-1)/2}}\b{(p-1)/2}{(p-1)/4} \e
(-1)^{\f{p-1}4}\cdot 2a\mod p.$$ Write
$P_{\f{p-1}2}(3)=(-1)^{\f{p-1}4}\cdot 2a+qp$. Then
$P_{\f{p-1}2}(3)^2\e 4a^2+(-1)^{\f{p-1}4}\cdot 4aqp\mod{p^2}$. By a
result due to Van Hamme [15], we have
$(\sum_{k=0}^{\f{p-1}2}\b{\f{p-1}2}k\b{\f{p-1}2+k}k)^2\e 4a^2-2p\mod
{p^2}.$ This together with (2.1) yields $P_{\f{p-1}2}(3)^2\e
4a^2-2p\mod{p^2}$. Hence $(-1)^{\f{p-1}4}\cdot 4aq\e -2\mod p$ and
so $P_{\f{p-1}2}(3)\e (-1)^{\f{p-1}4}\cdot 2a-\f
p{2(-1)^{(p-1)/4}a}\mod{p^2}.$ Now combining all the above we obtain
the result.
\par{\bf Remark 2.3} For a prime $p=4k+1=a^2+b^2$ with $a\e 1\mod
4$, the  congruence $\sum_{k=0}^{(p-1)/2}\b{2k}k^28^{-k}\e
 \sum_{k=0}^{(p-1)/2}\b{2k}k^2(-16)^{-k}
 \e  (-1)^{\f{p-1}4}\cdot 2a\mod p$
 was proved by Zhi-Wei Sun in [13], and he also conjectured
 $\sum_{k=0}^{(p-1)/2}\b{2k}k^28^{-k}\e
 \sum_{k=0}^{(p-1)/2}\b{2k}k^2(-16)^{-k} $ $
 \e  (-1)^{\f{p-1}4}(2a-\f p{2a})\mod {p^2}.$

 \pro{Theorem 2.10} Let $p$ be a prime of the
form $4k+1$ and so $p=a^2+b^2$ with $a,b\in\Bbb Z$ and $a\e 1\mod
4$. Then
$$\sum_{k=0}^{\f{p-1}4}\f{\b{4k}{2k}^2}{16^{2k}}
\e\f 12+ (-1)^{\f{p-1}4}a-(-1)^{\f{p-1}4}\f p{4a}\mod {p^2}.$$
\endpro
Proof. Since
$\sum_{k=0}^{\f{p-1}2}\b{2k}k^216^{-k}+\sum_{k=0}^{\f{p-1}2}
\b{2k}k^2
(-16)^{-k} =2\sum_{k=0}^{\f{p-1}4}\b{4k}{2k}^216^{-2k},$ by
(1.1) and Theorem 2.9 we deduce the result.
\par For a prime $p>3$ and $A,B,C\in \Bbb Z_p$ let
$\#E_p(y^2=x^3+Ax^2+Bx+C)$ be the number of points on the curve
$E_p:\ y^2=x^3+Ax^2+Bx+C$ over the field $\Bbb F_p$ of $p$ elements.
 \pro{Lemma 2.6 ([5])} Let $p>3$ be a prime and
$\lambda\in\Bbb Z_p$ with $\lambda\not\e 0,1\mod p$. Then
$$p+1-\#E_p(y^2=x(x-1)(x-\lambda))\e (-1)^{\f{p-1}2}
\sum_{k=0}^{\f{p-1}2}\b{\f{p-1}2+k}k\b{\f{p-1}2}k(-\lambda)^k\mod
p.$$
\endpro

\pro{Theorem 2.11} Let $p>3$ be a prime and $t\in\Bbb Z_p$. Then
$$P_{\f{p-1}2}(t)\e\sum_{k=0}^{\f{p-1}2}\b{2k}k^2\Ls{1-t}{32}^k
 \e -\Ls {-6}p\sum_{x=0}^{p-1}\Ls{x^3-3(t^2+3)x+2t(t^2-9)}p\mod
p.$$
\endpro
Proof. For $t\e \pm 1\mod p$ we have
$$\align &\sum_{x=0}^{p-1}\Ls{x^3-3(t^2+3)x+2t(t^2-9)}p\\&=
\sum_{x=0}^{p-1}\Ls{x^3-12x\mp
16}p=\sum_{x=0}^{p-1}\Ls{(2x)^3-12(2x)\mp 16}p =\Ls
2p\sum_{x=0}^{p-1}\Ls{x^3-3x\mp 2}p \\&=\Ls
2p\sum_{x=0}^{p-1}\Ls{(x\pm 1)^2(x\mp 2)}p =\Ls
2p\Big(\sum_{x=0}^{p-1}\Ls{x\mp 2}p-\Ls{\mp 3}p\Big) =-\Ls{\mp
6}p.\endalign$$ Thus applying (2.4) and the fact $P_n(\pm 1)=(\pm
1)^n$ (see [6]) we deduce the result.

\par Now assume $t\not\e \pm 1\mod p$. For $A,B,C\in\Bbb Z_p$, it is easily seen that (see for
example [12, pp.\;221-222])
$$\#E_p(y^2=x^3+Ax^2+Bx+C)=p+1+\sum_{x=0}^{p-1}\Ls{x^3+Ax^2+Bx+C}p.$$
Taking  $\lambda=(1-t)/2$ in Lemma 2.6 and applying the above and
(2.1) we see that
$$\align P_{\f{p-1}2}(t)&=\sum_{k=0}^{\f{p-1}2}\b{\f{p-1}2+k}k\b{\f{p-1}2}k\Ls{t-1}2^k
\\&\e (-1)^{\f{p-1}2}\big(p+1-\#E_p(y^2=x(x-1)(x-(1-t)/2))\big)\\&=
(-1)^{\f{p+1}2}\sum_{x=0}^{p-1}\Ls{x(x-1)(x-(1-t)/2)}p\mod
p.\endalign$$ Since
$$\align &\sum_{x=0}^{p-1}\Ls{x(x-1)(x-(1-t)/2)}p
\\&=\sum_{x=0}^{p-1}\Ls{\f x2(\f x2-1)(\f x2-\f{1-t}2)}p =\Ls 2p
\sum_{x=0}^{p-1}\Ls{x(x-2)(x+t-1)}p
\\&=\Ls 2p\sum_{x=0}^{p-1}\Ls{x^3+(t-3)x^2-2(t-1)x}p
\\&=\Ls
2p\sum_{x=0}^{p-1}\Ls{(x-\f{t-3}3)^3+(t-3)(x-\f{t-3}3)^2-2(t-1)(x-\f{t-3}3)}p
\\&=\Ls
2p\sum_{x=0}^{p-1}\Ls{x^3-\f{t^2+3}3x+\f{2t^3-18t}{27}}p =\Ls
2p\sum_{x=0}^{p-1}\Ls{(\f x3)^3-\f{t^2+3}3\cdot\f
x3+\f{2t^3-18t}{27}}p
\\&=\Ls
6p\sum_{x=0}^{p-1}\Ls{x^3-3(t^2+3)x+2t(t^2-9)}p,
\endalign$$
by the above and (2.4) we obtain the result. The proof is now
complete.
 \pro{Theorem 2.12} Let
$p>3$ be a prime.  Then
$$\align&P_{\f{p-1}2}(-31)\e\sum_{k=0}^{\f{p-1}2}\b{2k}k^2\e \sum_{k=0}^{\f{p-1}2}\f{\b{2k}k^2}{256^k}\e-\Ls
p3\sum_{x=0}^{p-1}\Ls{x^3-723x-7378}p\mod p,
\\&P_{\f{p-1}2}(33)\\&\e\sum_{k=0}^{\f{p-1}2}(-1)^k\b{2k}k^2\e(-1)^{\f{p-1}2}
\sum_{k=0}^{\f{p-1}2}\f{\b{2k}k^2}{(-256)^k} \e
(-1)^{\f{p+1}2}\sum_{x=0}^{p-1}\Ls{x^3-91x+330}p\mod p,
\\&P_{\f{p-1}2}(-15)\e\sum_{k=0}^{\f{p-1}2}\f{\b{2k}k^2}{2^k}\e\Ls 2p\sum_{k=0}^{\f{p-1}2}\f{\b{2k}k^2}{128^k}\e
(-1)^{\f{p+1}2}\sum_{x=0}^{p-1}\Ls{x^3-19x-30}p\mod p,
\\&P_{\f{p-1}2}(9)\e\sum_{k=0}^{\f{p-1}2}\f{\b{2k}k^2}{(-4)^k}\e(-1)^{\f{p-1}2}\sum_{k=0}^{\f{p-1}2}\f{\b{2k}k^2}{(-64)^k}\e
(-1)^{\f{p+1}2}\sum_{x=0}^{p-1}\Ls{x^3-7x+6}p\mod p,
\\&P_{\f{p-1}2}(5)\e\sum_{k=0}^{\f{p-1}2}\f{\b{2k}k^2}{(-8)^k}\e\Ls{-2}p
\sum_{k=0}^{\f{p-1}2}\f{\b{2k}k^2}{(-32)^k}\e -\Ls
p3\sum_{x=0}^{p-1}\Ls{x^3-21x+20}p\mod p.
\endalign$$
\endpro
Proof. Taking $t=-31$ in Theorem 2.11 and applying Theorem 2.7 we
see that
$$\align P_{\f{p-1}2}(-31)&\e\sum_{k=0}^{\f{p-1}2}\b{2k}k^2
\e \sum_{k=0}^{\f{p-1}2}\f{\b{2k}k^2}{256^k}\e-\Ls
{-6}p\sum_{x=0}^{p-1}\Ls{x^3-3(31^2+3)x-62(31^2-9)}p
\\&=-\Ls
{-6}p\sum_{x=0}^{p-1}\Ls{(2x)^3-4\cdot 723\cdot 2x-8\cdot 7378}p
\\&=-\Ls{-3}p \sum_{x=0}^{p-1}\Ls{x^3-723x-7378}p\mod p.
\endalign$$
Taking $t=33$ in Theorem 2.11 and applying Theorem 2.7 we see that
$$\align &P_{\f{p-1}2}(33)\\&\e\sum_{k=0}^{\f{p-1}2}(-1)^k\b{2k}k^2
\e \Ls{-16}p\sum_{k=0}^{\f{p-1}2}\f{\b{2k}k^2}{(-256)^k}\e-\Ls
{-6}p\sum_{x=0}^{p-1}\Ls{x^3-3276x+71280}p
\\&=-\Ls
{-6}p\sum_{x=0}^{p-1}\Ls{(6x)^3-36\cdot 91\cdot 6x+216\cdot 330}p
\\&=-\Ls{-1}p \sum_{x=0}^{p-1}\Ls{x^3-91x+330}p\mod p.
\endalign$$
The remaining congruences can be proved similarly.
\par For $a,b,n\in\Bbb N$, if $n=ax^2+by^2$ for some $x,y\in\Bbb Z$, we say that
$n=ax^2+by^2$. In 2003, Rodriguez-Villegas[11] posed many
conjectures on supercongruences. In particular, he conjectured that
for any prime $p>3$,
$$\sum_{k=0}^{p-1}\f{(4k)!}{256^kk!^4}\e\cases 4x^2-2p\mod {p^2}
&\t{if $p\e 1,3\mod 8$ and so $p=x^2+2y^2$,}
\\ 0\mod{p^2}&\t{if $p\e 5,7\mod 8$}
\endcases$$
and
$$\sum_{k=0}^{p-1}\f{\binom{2k}k^2\binom{3k}k}{108^k}
\e\cases 4x^2-2p\mod {p^2}&\t{if $3\mid p-1$ and so $p=x^2+3y^2$,}
\\0\mod {p^2}&\t{if $3\mid p-2$.}
\endcases$$
Recently the author's twin brother Zhi-Wei Sun ([13,14]) made a lot
of conjectures on supercongruences. In particular, he conjectured
that for a prime $p\not= 2,7$,
$$\sum_{k=0}^{p-1}\b{2k}k^3\e \sum_{k=0}^{p-1}\f{(4k)!}{81^kk!^4}\e\cases 4x^2-2p\mod {p^2}
&\t{if $\sls p7=1$ and so $p=x^2+7y^2$,}
\\ 0\mod{p^2}&\t{if $\sls p7=-1$.}
\endcases$$
 Inspired by their work, we pose  the
following conjectures.
 \pro{Conjecture 2.1}
Let $p>3$ be a prime. Then
$$\sum_{k=0}^{p-1}\f{(4k)!}{648^kk!^4}\e\cases 4x^2-2p\mod {p^2}
&\t{if $p\e 1\mod 4$ and $p=x^2+y^2$ with $2\nmid x$,}
\\ 0\mod{p^2}&\t{if $p\e 3\mod 4$.}
\endcases$$\endpro

\pro{Conjecture 2.2} Let $p>3$ be a prime. Then
$$\sum_{k=0}^{p-1}\f{(4k)!}{(-144)^kk!^4}\e\cases 4x^2-2p\mod {p^2}
&\t{if $p\e 1\mod 3$ and so $p=x^2+3y^2$,}
\\ 0\mod{p^2}&\t{if $p\e 2\mod 3$.}
\endcases$$\endpro

\pro{Conjecture 2.3} Let $p\not=2,3,7$ be a prime. Then
$$\sum_{k=0}^{p-1}\f{(4k)!}{(-3969)^kk!^4}\e\cases 4x^2-2p\mod {p^2}
&\t{if $p\e 1,2,4\mod 7$ and so $p=x^2+7y^2$,}
\\ 0\mod{p^2}&\t{if $p\e 3,5,6\mod 7$.}
\endcases$$\endpro

\pro{Conjecture 2.4} Let $p\not= 2,3,11$ be a prime. Then
$$\sum_{k=0}^{p-1}\f{(6k)!}{66^{3k}(3k)!k!^3}\e\cases \sls p{33}(4x^2-2p)\mod {p^2}
&\t{if $4\mid p-1$ and $p=x^2+y^2$ with $2\nmid x$,}
\\ 0\mod{p^2}&\t{if $4\mid p-3$.}
\endcases$$\endpro

\pro{Conjecture 2.5} Let $p>5$ be a prime. Then
$$\sum_{k=0}^{p-1}\f{(6k)!}{20^{3k}(3k)!k!^3}\e\cases \sls {-5}p(4x^2-2p)\mod {p^2}
&\t{if $p\e 1,3\mod 8$ and $p=x^2+2y^2$,}
\\ 0\mod{p^2}&\t{if $p\e 5,7\mod 8$.}
\endcases$$\endpro

\pro{Conjecture 2.6} Let $p>5$ be a prime. Then
$$\sum_{k=0}^{p-1}\f{(6k)!}{54000^k(3k)!k!^3}\e\cases \sls p5(4x^2-2p)\mod {p^2}
&\t{if $3\mid p-1$ and so $p=x^2+3y^2$,}
\\ 0\mod{p^2}&\t{if $3\mid p-2$.}
\endcases$$\endpro

\pro{Conjecture 2.7} Let $p>5$ be a prime. Then
$$\aligned&\sum_{k=0}^{p-1}\f{(6k)!}{(-12288000)^k(3k)!k!^3}\\&\e\cases \sls {10}p(L^2-2p)\mod {p^2}
&\t{if $p\e 1\mod 3$ and so $4p=L^2+27M^2$,}
\\ 0\mod{p^2}&\t{if $p\e 2\mod 3$.}
\endcases\endaligned$$\endpro

\pro{Conjecture 2.8} Let $p>7$ be a prime. Then
$$\sum_{k=0}^{p-1}\f{(6k)!}{(-15)^{3k}(3k)!k!^3}\e\cases \sls p{15}(4x^2-2p)\mod {p^2}
&\t{if $\sls p7=1$ and so $p=x^2+7y^2$,}
\\ 0\mod{p^2}&\t{if $\sls p7=-1$.}
\endcases$$\endpro

\pro{Conjecture 2.9} Let $p\not= 2,3,5,7,17$ be a prime. Then
$$\sum_{k=0}^{p-1}\f{(6k)!}{255^{3k}(3k)!k!^3}\e\cases \sls p{255}(4x^2-2p)\mod {p^2}
&\t{if $\sls p7=1$ and so $p=x^2+7y^2$,}
\\ 0\mod{p^2}&\t{if $\sls p7=-1$.}
\endcases$$\endpro

\pro{Conjecture 2.10} Let $p>3$ be a prime. Then
$$\sum_{k=0}^{p-1}\f{\binom{2k}k^2\binom{3k}k}{1458^k}\e\cases 4x^2-2p\mod {p^2}
&\t{if $p\e 1\mod 3$ and so $p=x^2+3y^2$,}
\\ 0\mod{p^2}&\t{if $p\e 2\mod 3$.}
\endcases$$\endpro

\pro{Conjecture 2.11} Let $p>5$ be a prime. Then
$$\sum_{k=0}^{p-1}\f{\binom{2k}k^2\binom{3k}k}{15^{3k}}\e\cases 4x^2-2p\mod {p^2}
&\t{if $p\e 1,4\mod {15}$ and so $p=x^2+15y^2$,}
\\2p-12x^2\mod {p^2}
&\t{if $p\e 2,8\mod {15}$ and so $p=3x^2+5y^2$,}
\\ 0\mod{p^2}&\t{if $p\e 7,11,13,14\mod {15}$.}
\endcases$$\endpro

\pro{Conjecture 2.12} Let $p>5$ be a prime. Then
$$\sum_{k=0}^{p-1}\f{\binom{2k}k^2\binom{3k}k}{(-8640)^k}\e
\cases 4x^2-2p\mod{p^2}&\t{if $3\mid p-1$, $p=x^2+3y^2$ and $5\mid
xy$,}
\\p-2x^2\pm 6xy\mod{p^2}&\t{if $3\mid p-1$, $p=x^2+3y^2$, $5\nmid xy$}
\\&\qq\t{and $x\e \pm y,\pm 2y\mod 5$,}
\\0\mod{p^2}&\t{if $3\mid p-2$.}\endcases$$
\endpro

\pro{Conjecture 2.13} Let $p>3$ be a prime. Then
$$\sum_{k=0}^{p-1}\f{(3k)!}{54^k\cdot k!^3}
\e\cases \sls x3(2x-\f p{2x})\mod {p^2} &\t{if $3\mid p-1$ and so
$p=x^2+3y^2$,}
\\ 0\mod {p^2}&\t{if $3\mid p-2$.}
\endcases$$\endpro

\enddocument
\bye